\documentclass[12 pt]{amsart}
\usepackage{amscd,amsfonts,amssymb,amsmath}
\usepackage{hyperref}
\usepackage[cp1251]{inputenc}
\usepackage[T2A]{fontenc}

\usepackage{epsfig}
\newtheorem{theorem}{Theorem}[section]
\newtheorem{corollary}[theorem]{Corollary}

\newtheorem{proposition}[theorem]{Proposition}
\theoremstyle{definition}
\newtheorem{conjecture}[theorem]{Conjecture}
\newtheorem{question}[theorem]{Question}

\newtheorem{definition}[theorem]{Definition}
\newtheorem{example}[theorem]{Example}

\numberwithin{equation}{subsection}
\newtheorem*{ack}{Acknowledgement}
\usepackage[all,cmtip]{xy}
\usepackage{xcolor}
\usepackage{xfrac}
\usepackage{pb-diagram}

\def\bea{\begin{eqnarray}}
\def\eea{\end{eqnarray}}

\begin{document}
\title{Multivalued groups and  Newton polyhedron}
\author{Valeriy G. Bardakov}
\author{Tatyana A. Kozlovskaya}

\date{\today}
\address{Sobolev Institute of Mathematics and Novosibirsk State University, Novosibirsk 630090, Russia.}
\address{Novosibirsk State Agrarian University, Dobrolyubova street, 160, Novosibirsk, 630039, Russia.}
\address{Regional Scientific and Educational Mathematical Center of Tomsk State University,
36 Lenin Ave., 14, 634050, Tomsk, Russia.}
\email{bardakov@math.nsc.ru}

\address{Regional Scientific and Educational Mathematical Center of Tomsk State University, 36 Lenin Ave., 14, 634050, Tomsk, Russia.}
\email{t.kozlovskaya@math.tsu.ru}


\begin{abstract}
On the set of complex number $\mathbb{C}$ it is possible to define $n$-valued group for any positive integer $n$. The $n$-multiplication defines a symmetric polynomial $p_n = p_n(x, y, z)$ with integer coefficients. By the theorem on symmetric polynomials, one can present $p_n$ as polynomial in elementary symmetric polynomials $e_1$, $e_2$, $e_3$. V.~M.~Buchstaber formulated a question on description coefficients of this polynomial. Also, he formulated the next question:  How to describe the Newton polyhedron of $p_n$?  In the present paper we find all coefficients of $p_n$ under monomials of the form $e_1^i e_2^j$ and prove that  the Newton polyhedron of $p_n$ is an   right  triangle.
\end{abstract}
\maketitle

\section{Introduction}
One  branch  of Abstract Algebra is studying algebraic systems with multivalued operations. Solutions of the Yang-Baxter equation (2-simplex equation) and its generalization $n$-simplex equations, $n \geq 3$, are examples of multivalued operations. 
In 1971, V.~M.~Buchstaber and S.~P.~Novikov \cite{BN} introduced a construction, suggested by
the theory of characteristic classes of vector bundles, in which the product of each pair of elements is an $n$-multi-set, the set of $n$ points with multiplicities. This construction leads to the notion of $n$-valued group.

A good survey on $n$-valued groups and its applications  can be found in \cite{B}. In Section 5 of this paper  on the set of complex numbers $\mathbb{C}$ were constructed $n$-valued groups   for any natural $n$. The $n$-valued multiplication is describe by the polynomials $p_n = p_n(z;x,y)$ which are 
 $x, y, z$--symmetric polynomials with integer coefficients. If we introduce elementary symmetric polynomials
 $$
e_1 = x + y + z, ~~e_2 = xy + yz + zx,~~ e_3 = xyz,
$$
then $p_n = P_n(e_1, e_2, e_3)$ is a polynomial with integer coefficients in variables $e_1$, $e_2$, $e_3$. In \cite{B} were  formulated two questions on the description the coefficients of $P_n$ and a question on the Newton polyhedron of $p_n$.

In the present paper we find all coefficients at $e_1^{i} e_2^{j}$ in $P_n$. It gives particular answer on the first two questions. 
Also, we  prove that if  $f = f(x_1, x_2, \ldots, x_n) \in \mathbb{Z}[x_1, x_2, \ldots, x_n]$ is a symmetric homogeneous polynomial of degree $k$, which  contains a monomial $a x_1^k$ for some non-zero $a$, then  its Newton polyhedron  is the $k \Delta^{n-1}$-simplex.
From this theorem follows 
 that  the Newton polyhedron of $p_n$ is the right triangle with side which depend on $n$. It is the full answer on the third question.
 
 At the end of the paper we formulate some open questions.

\bigskip


\section{Multivalued groups and Buchstaber's questions}

\subsection{Multivalued groups}
Recall definitions and some facts from the theory of multivalued groups (see, for example, \cite{B}).

Let $X$ be a non-empty set. An $n$-{\it valued multiplication} on $X$ is a map
$$
\mu \colon X \times X \to (X)^n = Sym^nX,~~\mu(x, y) = x * y = [z_1, z_2, \ldots, z_n],~~z_k = (x*y)_k,
$$
where $(X)^n = Sym^nX$ is the $n$-th symmetric power of $X$, that is the quotient $X^n / S_n$ of the Cartesian power $X^n$ under the action of $S_n$ by permutations of components. 
The next axioms are  natural generalizations of the classical axioms of group multiplication.

{\it Associativity}. The $n^2$-multi-sets:
$$
[x * (y * z)_1, x * (y * z)_2,  \ldots, x * (y * z)_n], ~~[ (x * y)_1 * z, (x * y)_2 * z,  \ldots, (x * y)_n * z]
$$
are equal for all $x, y, z \in X$.

{\it Unit}. An element $e \in X$ such that 
$$
e * x = x * e = [x, x, \ldots, x]
$$
 for all $x \in X$.

{\it Inverse}. A map $inv \colon X \to X$ such that 
$$
e \in inv(x) * x~\mbox{and}~e \in x * inv(x)
$$
 for all $x \in X$.

The map $\mu$ defines $n$-{\it valued group structure} $\mathcal{X} = (X, \mu, e, inv)$ on $X$ if it is associative, has a unit and an inverse.

\medskip

 Let $\mu$ be the multiplication
$$
\mu \colon \mathbb{C} \times \mathbb{C} \to (\mathbb{C})^n
$$
that is defined by the formula
$$
\mu(x, y) =  x * y = [ (\sqrt[n]{x} + \epsilon^r \sqrt[n]{y})^n,~~1 \leq r \leq n],
$$
where $\epsilon$ is a primitive $n$-th root of unity.
This multiplication endows $\mathbb{C}$ with the structure of an $n$-valued group with the unit $e = 0$. The inverse element is given by the map $inv(x) = (-1)^n x$.

The $n$-valued multiplication is described by the polynomials
$$
p_n = p_n(z; x, y) = \prod^n_{k=1} (z - (inv(x) * inv(y))_k),
$$
whence the product $x * y$ is defined by $z$-roots of the equation $p_n = 0$. The polynomials 
$p_n(z; x, y)$ are $x, y, z$--symmetric polynomials with integral coefficients, e.g.,
$$
p_1 = x + y + z, ~~p_2 = (x + y + z)^2 - 4(xy + yz + zx).
$$
Set
$$
e_1 = x + y + z, ~~e_2 = xy + yz + zx,~~ e_3 = xyz.
$$
Then

$p_1 = e_1,$\\

$p_2 = e_1^2 - 2^2 \,  e_2,$\\

$p_3 = e_1^3 - 3^3  \,  e_3,$\\

$p_4 = e_1^4 - 2^3  \,  e_1^2 e_2 + 2^4  \,  e_2^2 - 2^7  \,  e_1 e_3,$\\

$p_5 = e_1^5 - 5^4  \,  e_1^2 e_3 + 5^5  \,  e_2 e_3,$\\

$p_6 = e_1^6 - 2^2 \cdot 3  \,  e_1^4 e_2 + 2^4 \cdot 3  \,  e_1^2 e_2^2 - 2^6  \,  e_2^3 -   2 \cdot 3^4 \cdot 17  \,  e_1^3 e_3  - 2^3  \cdot 3^4 \cdot 19   \,  e_1 e_2 e_3 + 3^3 \cdot 19^3  \,  e_3^2 ,$\\

$p_7 = e_1^7 -  5 \cdot 7^4  \,  e_1^4 e_3 +    2 \cdot 7^6  \,   e_1^2 e_2 e_3  - 7^7    \,  e_2^2 e_3 + 7^8  \,  e_1 e_3^2.$\\

The next questions were formulated in \cite{B}.

\begin{question} \label{BQ}
(1) What is the relationship between prime factors of $n$ and prime factors of the
coefficients of the polynomials $p_n$?

(2) How to distinguish the monomials that have zero coefficient? 

(3) How to describe the Newton polyhedron of $p_n$?
\end{question}

\bigskip


\section{Coefficients  and the Newton polyhedron  of $p_n$} 

Since $p_n$ is a symmetric homogeneous  polynomial of degree $n$, by the theorem on symmetric polynomials we can present $p_n$ as a polynomial on the elementary symmetric polynomials $e_1$, $e_2$, and $e_3$,
$$
p_n = \sum_{\substack{k_1 \geq k_2 \geq k_3 \geq 0\\ k_1 + k_2 + k_3 = n}} A_{k_1,k_2,k_3} e_1^{k_1-k_2} e_2^{k_2-k_3} e_3^{k_3} \in \mathbb{Z}[e_1, e_2, e_3].
$$
The main problem is to find the coefficients $A_{k_1,k_2,k_3}$. 

We can write $p_n$ in the form
$$
p_n =  \prod^n_{k=1} \left(z - \left((inv(x) * inv(y)\right)_k\right) = \prod^n_{k=1} \left(z - ((-1)^n x * (-1)^n y)_k\right) =
$$
$$
= \prod^n_{k=1} \left(z - \left(\sqrt[n]{(-1)^n x} +  \epsilon^k \sqrt[n]{(-1)^n y}\right)^n\right).
$$
If $y = 0$, then
$$
\bar{p}_n = p_n(z; x, 0) = \prod^n_{k=1} \left(z - (\sqrt[n]{(-1)^n x})^n\right) = \prod^n_{k=1} \left(z - (-1)^n x\right) =   \left(z - (-1)^n x\right)^n.
$$

Denote by
$$
\bar{e}_1 = e_1(z;x,0) = x + z,~~\bar{e}_2 = e_2(z;x,0) = zx.
$$
We see that $e_3(z;x,0) = 0$.

The next proposition gives particular answers on the first two questions.

\begin{proposition}
1) If $n$ is odd, then all $A_{k_1,k_2,0}$, $k_2 \not= 0$, are zero, i.e. in this case $p_n$ does not contains monomials $e_1^i e_2^j$, $j > 0$.

2) If $n = 2k$ is even, then the coefficient $A_{2k-i,i,0}$ at $e_1^{2(k-i)} e_2^i$, is equal to
$$
A_{2k-i,i,0} = (-4)^i C_k^i =  (-4)^i \frac{k!}{i! (k-i)!},~~~i = 1, 2, \ldots, k.
$$
\end{proposition}

\begin{proof}
1) If $n$ is odd, then 
$$
\bar{p}_n =  (z + x)^n = \bar{e}_1^n.
$$
It means that in $p_n$ all coefficients  $A_{k_1,k_2,0}$, where $k_1 \geq k_2 >0$ and $k_1 + k_2 = n$ are zero.

2) If $n=2k$ is even, then
$$
\bar{p}_n =  (z - x)^n = (\bar{e}_1^2 - 4 \bar{e}_2)^k = \sum_{i=0}^k (-4)^i C_k^i (\bar{e}_1^2)^{k-i} (\bar{e}_2)^i.
$$
Hence, in $p_n$ we have found the next coefficients
$$
A_{2k-i,i,0} = (-4)^i C_k^i,~~~i = 1, 2, \ldots, k.
$$
\end{proof}

\begin{example}
From this proposition follows that for even $n$ hold\\

$\bar{p}_2 =   \bar{e}_1^2 - 2^2   \,  \bar{e}_2,$\\

$\bar{p}_4 =   \bar{e}_1^4 - 2^3   \,  \bar{e}_1^2 \bar{e}_2 + 2^4   \,  \bar{e}_2^2,$\\

$\bar{p}_6 =   \bar{e}_1^6 - 2^2 \cdot 3   \,  \bar{e}_1^4 \bar{e}_2 +   2^4 \cdot 3   \,  \bar{e}_1^2 \bar{e}_2^2 - 2^6   \,  \bar{e}_2^3,$\\

$\bar{p}_8 =   \bar{e}_1^8 - 2^4   \,  \bar{e}_1^6 \bar{e}_2 +   2^5 \cdot 3   \,  \bar{e}_1^4 \bar{e}_2^2 - 2^8   \,  \bar{e}_1^2 \bar{e}_2^3  + 2^8   \,  \bar{e}_2^4.$
\end{example}

It is easy to see that for even $n$ all coefficients of $\bar{p}_n$ except the coefficient at $\bar{e}_1^n$ are even. It is not true for polynomials $p_n$, as example $p_6$ shows. We can formulate

\begin{conjecture}
1) If $n = p^m$ is a power of a prime $p$, then all coefficients, except the coefficient at $e_1^n$ are divided into $p$. 

2) If $n$ is even, then all coefficients $A_{k_1,k_2,k_3}$ are non-zero. 
\end{conjecture}

\subsection{Newton polyhedron} In this subsection we give a full answer on the third question in  \ref{BQ}.
Recall the need definition. Let 
$$
f = f(x_1, x_2, \ldots, x_n) = \sum a_{i_1 \ldots i_n}  x_1^{i_1} \ldots x_n^{i_n} \in \mathbb{Z}[x_1, x_2, \ldots, x_n]
$$
 be a polynomial with integer coefficients. Denote by $I_f$ the set of multi indexes $(i_1, \ldots, i_n)$ such that $a_{i_1 \ldots i_n} \not= 0$. The convex hull 
$$
N_f = Conv(I_f) \subset \mathbb{R}^n
$$
is said to be a {\it Newton polyhedron} of $f$.

To find the  Newton polyhedrons for the polynomials $p_n$, consider them for small $n$,

$p_1 = x + y + z,$\\

$p_2 = x^2 + y^2 + z^2 - 2 x y - 2 y z - 2 z x,$\\

$p_3 =(z + x + y)^3 - 27 x y z,$\\

$p_4 = ((x + y + z)^2 - 4(x y + y z + z x))^2 - 2^7(x + y + z)x y z = p_2^2 - 2^7 p_1 x y z.$\\

Denote by $N_i \subset \mathbb{R}^3$ the Newton polyhedron for $p_i$. Then 

-- $N_1$ is the right triangle $A_1 B_1 C_1$ with the vertices $A_1 = (1, 0, 0)$, $B_1 = (0, 1, 0)$,
$C_1 = (0, 0, 1)$;

-- $N_2$ is the right triangle $A_2 B_2 C_2$ with the vertices $A_2 = (2, 0, 0)$, $B_2 = (0, 2, 0)$,
$C_2 = (0, 0, 2)$;

-- $N_3$ is the right triangle $A_3 B_3 C_3$ with the vertices $A_3 = (3, 0, 0)$, $B_3 = (0, 3, 0)$,
$C_3 = (0, 0, 3)$;

-- $N_4$ is the right triangle $A_4 B_4 C_4$ with the vertices $A_4 = (4, 0, 0)$, $B_4 = (0, 4, 0)$,
$C_4 = (0, 0, 4)$.

To describe $N_k$ for $k >2$ we introduce the next definition.

\begin{definition}
Let $k$ be a positive integer. The {\it standard $n$-simplex of size} $k$ is the subset of $\mathbb{R}^{n+1}$ given by
$$
k \Delta^n = \left\{ (t_0, t_1, \ldots, t_n) \in \mathbb{R}^{n+1}~|~\sum_{i=0}^n t_i = k ~\mbox{and}~t_i \geq 0~ \mbox{for}~i = 0, 1, \ldots, n \right\}.
$$
For simplicity we shall call the standard $n$-simplex of size $k$ by $k \Delta^n$-simplex.
\end{definition}
For $k=1$ we get the definition of the standard $n$-simplex (or unit simplex).

The $k \Delta^n$-simplex has $n+1$ vertices,
$$
E_0 =  (k, 0, 0, \ldots, 0, 0), E_1 =  (0, k, 0,  \ldots, 0, 0), \ldots E_n =  (0, 0, 0, \ldots, 0, k).
$$

Now we are ready to prove the main result of the present subsection.

\begin{theorem} 
 Let $f = f(x_1, x_2, \ldots, x_n) \in \mathbb{Z}[x_1, x_2, \ldots, x_n]$ be a symmetric homogeneous polynomial of degree $k$, which  contains a monomial $a x_1^k$ for some non-zero $a$. Then  its Newton polyhedron $N_f$ is the $k \Delta^{n-1}$-simplex.
\end{theorem}

\begin{proof}
Since $a x_1^k$ is a monomial of $f$ and $f$ is symmetric, it contains monomials $a x_i^k$ for all $i = 1, 2, \ldots, n$. Hence,   $N_f$ contains the vertices 
$$
E_0 =  (k, 0, 0, \ldots, 0, 0), E_1 =  (0, k, 0,  \ldots, 0, 0), \ldots E_{n-1} =  (0, 0, 0, \ldots, 0, k) \in \mathbb{R}^n
$$
and hence contains  $k \Delta^{n-1}$-simplex. Let us show that any other vertex of $N_f$, which corresponds a monomial in $f$ lies in this simplex. Indeed, any such monomial has the form 
$$
b x_1^{k_1} x_2^{k_2} \ldots x_n^{k_n},~~b \in \mathbb{R},~~ b \not= 0.
$$
Since 
$$
k_1 + k_2 + \ldots + k_n = k,~~k_i \geq 0~ \mbox{for}~i = 1, 2, \ldots, n,
$$
the corresponding vertex lies in  $k \Delta^{n-1}$-simplex.
\end{proof}

We seen that the polynomial $p_k$ is homogeneous and has the form $p_k = e_1^k + \ldots$. Hence, from the theorem follows an answer on the third question of V. M. Buchstaber.

\begin{corollary}
The Newton polyhedron that corresponds to the polynomial $p_k(x, y, z)$, $k \geq 1$, is the $k \Delta^{2}$-simplex that is a right triangle  with sides of length $\sqrt[]{2} \, k$.
\end{corollary}

\bigskip


\section{Some open questions}

The following questions seem interesting.

\begin{question}
\begin{enumerate}
\item Let $f = f(x_1, x_2, \ldots, x_n) \in \mathbb{Z}[x_1, x_2, \ldots, x_n]$ be a symmetric polynomial, $N_f$ is its Newton polyhedron. Let us present $f$ as a polynomial on elementary symmetric polynomial, $f =F[e_1, e_2, \ldots, e_n] \in \mathbb{Z}[e_1, e_2, \ldots, e_n]$ and construct  its Newton polyhedron $N_F$. What is connection between $N_f$ and $N_F$?

\item Let $f = f(x_1, x_2, \ldots, x_n) \in \mathbb{Z}[x_1, x_2, \ldots, x_n]$ be a symmetric polynomial of degree $k$, which does not contain $x_1^k$. What can we say on its Newton polyhedron $N_f$?

\item Since there exists homomorphism of one multivalued group to another multivalued group and define a kernel of this homomorphism, we can say on extensions of multivalued groups. Construct a theory of extensions multivalued groups. 

\item Is it possible to define (co)homology for  multivalued groups?

\end{enumerate}
\end{question}

\begin{ack}
This paper is supported by the Ministry of Science and Higher Education of Russia (agreement  No.  075-02-2022-884). The authors thank  V. M. Buchstaber and D. V. Talalaev for interesting discussions and useful suggestions.
\end{ack}

\end{document}